\documentclass[centertags,leqno]{article}

\usepackage{amsmath}
\usepackage{amssymb}
\usepackage{latexsym}
\usepackage{amsxtra}
\usepackage{amscd}
\usepackage{theorem}

\usepackage{curves}
\usepackage{epic}
\usepackage{eepic}

\theoremstyle{plain}

{\theorembodyfont{\slshape}

        \newtheorem{thm}{Theorem}[section]

        \newtheorem{prop}[thm]{Proposition}
        
}

{\theorembodyfont{\rmfamily}

        \newtheorem{defn}[thm]{Definition}

}

\renewcommand{\em}{\sl}

\makeatletter
\renewcommand{\subsection}{\@startsection{subsection}{2}%
        {\z@}{-3.25ex plus -1ex minus-.2ex}{-1em}{\bf}}
\makeatother

\newcommand{\PP}{\mathbb{P}}
\newcommand{\HH}{\mathbb{H}}
\newcommand{\ZZ}{\mathbb{Z}}
\newcommand{\CC}{\mathbb{C}}

\newcommand{\QQ}{\mathbb{Q}}

\newcommand{\FF}{\mathbb{F}}

\newcommand{\X}{\mathcal{X}}
\newcommand{\Y}{\mathcal{Y}}

\newcommand{\OO}{\mathcal{O}}

\newcommand{\p}{\mathfrak{p}}

\newcommand{\SL}{{\rm SL}}

\newcommand{\PSL}{{\rm PSL}}
\newcommand{\Gal}{{\rm Gal}}
\newcommand{\Aut}{{\rm Aut}}

\newcommand{\Ind}{{\rm Ind}}

\newcommand{\Kb}{\bar{K}}
\newcommand{\QQb}{\bar{\QQ}}
\newcommand{\Xb}{\bar{X}}
\newcommand{\Yb}{\bar{Y}}
\newcommand{\Zb}{\bar{Z}}
\newcommand{\fb}{\bar{f}}

\newcommand{\an}{^{\rm an}}
\newcommand{\hor}{^{\rm hor}}

\newcommand{\aux}{^{\rm aux}}

\newcommand{\stab}{^{\rm\scriptscriptstyle st}}

\newcommand{\inj}{\hookrightarrow}
\newcommand{\To}{\;\longrightarrow\;}

\title{Stable reduction of three point covers}

\author{
   Stefan Wewers
 \\ Universit\"at Bonn}

\date{}

\begin{document}

\maketitle

This note gives a survey on some results related to the stable reduction of
three point covers, which were the topic of my talk at the {\em Journ\'ees
  Arithmetiques} 2003 in Graz. I thank the organizers of this very nice
conference for inviting me and letting me give a talk.

\section{Three point covers} \label{threepoint}

\subsection{Ramification in the field of moduli} \label{fom}

Let $X$ be a smooth projective curve over $\CC$. A celebrated theorem
of Belyi states that $X$ can be defined over a number field $K$ if and
only if there exists a rational function $f$ on $X$ with exactly three
critical values, see \cite{Belyi}, \cite{KoeckBelyi}. If such a
function $f$ exists, we can normalize it in such a way that the
critical values are $0$, $1$ and $\infty$. After this normalization,
we may view $f$ as a finite cover $f:X\to\PP^1$ which is \'etale over
$\PP^1-\{0,1,\infty\}$. We call $f$ a {\em three point cover}.
Another common name for $f$ is {\em Belyi map}. The  {\em monodromy
  group} of $f$ is defined as the Galois group of the Galois closure
of $f$.  

Let $f:X\to\PP^1$ be a three point cover. By the `obvious direction'
of Belyi's theorem, $f$ can be defined over the field $\QQb$ of
algebraic numbers. Therefore, for every element
$\sigma\in\Gal(\QQb/\QQ)$ of the absolute Galois group of $\QQ$ we
obtain a conjugate three point cover
$f^\sigma:X^\sigma\to\PP^1$, which may or may not be isomorphic
to $f$. This yields a continuous action of $\Gal(\QQb/\QQ)$ on the set
of isomorphism classes of three point covers. Hence we can
associate to $f$ the number field $K$ such that $\Gal(\QQb/K)$ is
precisely the stabilizer of the isomorphism class of $f$. The field
$K$ is called the {\em field of moduli} of $f$. Under certain extra
assumptions on $f$, the field $K$ is the smallest field of definition
of the cover $f$, see \cite{DebDou1}

Three point covers are determined, up to isomorphism, by finite,
purely combinatorial data -- e.g.\ by a {\em dessin d'enfants}
\cite{SchnepsDessins}.  It is an interesting problem to describe the
field of moduli of a three point cover in terms of these data. There
are only few results of a general nature on this problem. The aim of
this note is to explain certain results leading to the following
theorem, proved in \cite{bad}.

\begin{thm} \label{tamethm}
  Let $f:X\to\PP^1$ be a three point cover, with field of moduli $K$
  and monodromy group $G$. Let $p$ be a prime number such that $p^2$
  does not divide the order of $G$. Then $p$ is at most tamely
  ramified in the extension $K/\QQ$.
\end{thm}

If the prime $p$ does not divide the order of $G$ then $p$ is even
unramified in the extension $K/\QQ$, by a well known theorem of
Beckmann \cite{Beckmann89}. Both Beckmann's result and Theorem
\ref{tamethm} rely on an analysis of the reduction of $f$ at the prime
ideals $\p$ of $K$ dividing $p$. The results leading to Theorem
\ref{tamethm} were mainly inspired by Raynaud's paper
\cite{Raynaud98}.

\subsection{Good reduction} \label{good}

For a discussion of the relation between the ramification in the
extension $K/\QQ$ and the reduction behaviour of $f$, it is convenient
to localize at a prime ideal $\p$ of $K$ dividing $p$. In other words,
we may replace the extension $K/\QQ$ by a finite extension of $p$-adic
fields.

Fix a prime number $p$ and let $K_0$ denote the completion of the
maximal unramified extension of $\QQ_p$.  From now on, the letter $K$
will always denote a finite extension of $K_0$.  Note that $K$ is
complete with respect to a discrete valuation $v$, with residue field
$k=\bar{\FF}_p$. 

Let $f:X\to\PP^1_K$ be a three point cover, defined over a finite
extension $K/K_0$. Let $G$ denote the monodromy group of $f$. The
theorem of Beckmann mentioned above can be reduced to the following
statement: if $p$ does not divide the order of $G$ then $f$ can be
defined over $K_0$. This latter statement is, more or less, a direct
consequence of Grothendieck's theory of the tame fundamental group.
Indeed, if $p$ does not divide the order of $G$ then $f$ extends, by
Grothendieck's theory, to a tame cover $f_{\OO_K}:\X\to\PP^1_{\OO_K}$,
ramified only along the sections $0$, $1$, $\infty$. Let
$\fb:\Xb\to\PP^1_k$ denote the special fibre of this map; it is a tame
cover, ramified at most at $0$, $1$, $\infty$.  By the deformation
theory of tame covers, there exists a tame cover
$f_{\OO_{K_0}}:\X_0\to\PP^1_{\OO_{K_0}}$, ramified at most along $0$,
$1$, $\infty$ lifting $\fb$. Moreover, such a lift is unique. It
follows that the generic fibre of $f_{\OO_{K_0}}$ is a model of $f$
over $K_0$.

To see how Beckmann's Theorem follows from the above, let
$f:X\to\PP^1$ be a three point cover, defined over $\QQb$. Let
$\p$ be a place of $\QQb$ whose residue characteristic is prime to the
order of the monodromy group of $f$. Let $\sigma\in\Gal(\QQb/\QQ)$ be
an element of the inertia group of $\p$. We claim that $f^\sigma\cong
f$.  Clearly, this claim implies Beckmann's Theorem.

To prove the claim, let $K_0$ be as above. The place $\p$ gives rise
to an embedding $\QQb\inj\Kb_0$. Moreover, there exists a (unique)
element $\tau\in\Gal(\Kb_0/K_0)$ with $\tau|_{\QQb}=\sigma$. It
follows from the good reduction result discussed above that the three point
cover $f_{\Kb_0}:=f\otimes\Kb_0$ can be defined over $K_0$. Hence we
have $f_{\Kb_0}^\tau\cong f_{\Kb_0}$, which implies $f^\sigma\cong f$.
This proves the claim.

\section{Stable reduction} \label{reduction}

\subsection{The stable model} \label{stable}

Let $K_0$ be as in \S \ref{good}, and let $f:X\to\PP^1_K$ be a three point
cover, defined over a finite extension $K/K_0$.  If $p$ divides the
order of the monodromy group, then $f$ may have bad reduction and it
may not be possible to define $f$ over $K_0$. For the purpose of
studying this situation, it is no restriction to make the following
additional assumptions.
\begin{itemize}
\item
  The cover $f:X\to\PP^1_K$ is Galois, with Galois group $G$ (replace
  $f$ by its Galois closure).
\item
  The curve $X$ has semistable reduction (replace $K$ by a finite
  extension). 
\end{itemize}
In the second point, we have used the Semistable Reduction Theorem, see e.g.\ 
\cite{DeligneMumford}. For simplicity, we shall also assume
that the genus of $X$ is at least two. (Three point Galois covers of genus
$\leq 1$ can be classified and treated separately.)  Let $\X$ denote
the stable model of $X$, i.e.\ the minimal semistable model of $X$ over the
ring of integers of $K$, see \cite{DeligneMumford}. By uniqueness of the
stable model, the action of $G$ on $X$ extends to $\X$. Let $\Y:=\X/G$ be the
quotient scheme. It is shown in \cite{Raynaudfest}, Appendice, that $\Y$ is
again a semistable curve over $\OO_K$.

\begin{defn} \label{stabdef}
  The morphism $f\stab:\X\to\Y$ is called the {\em stable model} of
  the three point cover $f$. Its special fibre $\fb:\Xb\to\Yb$ is
  called the {\em stable reduction} of $f$. If $\fb$ is a separable
  and tamely ramified map between smooth curves, then we say that $f$
  has {\em good reduction}. Otherwise, we say that $f$ has bad
  reduction.
\end{defn}

Initiated by a series of papers by Raynaud \cite{Raynaudfest},
\cite{Raynaud94}, \cite{Raynaud98}, several authors have studied the
stable reduction of covers of curves (the case of three point covers
is just a special case).  For an overview of their results and a more
extensive list of references, see \cite{LiuOverview}. In this note, we
shall focus on the results of \cite{bad}, and on results which
inspired this work (mainly \cite{Raynaud98}, \cite{YannickArbres},
\cite{special}).

\subsection{Bad reduction} \label{badreduction}

Let $\fb:\Xb\to\Yb$ be the stable reduction of a three point cover
$f:X\to\PP^1_K$.  Let $(\Yb_i)$ be the list of all irreducible
components of the curve $\Yb$. Since the generic fibre of $\Y$ is just
the projective line, the components $\Yb_i$ are all smooth curves of
genus $0$. Moreover, the graph of components of $\Yb$ (whose vertices
are the components $\Yb_i$ and whose edges are the singular points) is
a tree.  For each index $i$, we fix an irreducible component $\Xb_i$ of
$\Xb$ such that $\fb(\Xb_i)=\Yb_i$. Let $\fb_i:\Xb_i\to\Yb_i$ denote
the restriction of $\fb$ to $\Xb_i$. Let $G_i\subset G$ denote the
stabilizer of the component $\Xb_i$.

The component $\Yb_i$ corresponds to a discrete valuation $v_i$ of the
function field $K(Y)$ of $Y=\PP^1_K$ whose residue field is the
function field of $\Yb_i$. The choice of $\Xb_i$ corresponds to the
choice of a valuation $w_i$ of the function field $K(X)$ of $X$
extending $v_i$, and the map $\fb_i$ corresponds to the residue field
extension of $w_i|v_i$. The group $G_i$ is simply the decomposition
group of $w_i$ in the Galois extension $K(X)/K(Y)$. Let $I_i\lhd G_i$
denote the corresponding inertia group.

By construction, the curve $\Xb$ is reduced. It follows that the
ramification index of the extension of valuations $e(w_i/v_i)$ is equal
to one. This implies that the inertia group $I_i$ is a $p$-group whose
order is equal to the degree of inseparability of the extension of
residue fields.

We say that $\Yb_i$ is a {\em good component} if the map $\fb_i$ is
separable. By what we have said above, this holds if and only if
$I_i=1$, i.e.\ the valuation $v_i$ is unramified in the 
extension $K(X)/K(Y)$. If this is the case, then $\fb_i:\Xb_i\to\Yb_i$
is a Galois cover with Galois group $G_i$.

If $\fb_i$ is not separable we say that $\Yb_i$ is a {\em bad
  component}. The map $\fb_i$ factors as the composition of a purely
inseparable map $\Xb_i\to\Zb_i$ of degree $|I_i|$ and a Galois cover
$\Zb_i\to\Xb_i$ with Galois group $G_i/I_i$.

Note that $K(Y)=K(t)$, where $t$ is the standard parameter on $\PP^1$.
To simplify the exposition, we shall make the following additional
assumption: there is a (necessarily unique) component $\Yb_0$ of $\Yb$
which corresponds to the Gauss valuation on $K(t)$ with respect to the
parameter $t$. This component is called the {\em original component}.
It is canonically isomorphic to $\PP^1_k$. (By making this assumption,
we exclude the case where the cover $f$ has bad reduction but the
curve $X$ has good reduction. In \cite{bad}, Definition 2.1, this is
called the {\em exceptional case}.) Let $\Yb_1,\ldots,\Yb_r$ be the
components of $\Yb$ different from $\Yb_0$.

\begin{figure}
\begin{center}
\unitlength3mm

\begin{picture}(30,16)

\put(6,4){\line(1,0){16}}
\put(6,10){\line(1,0){16}}
\put(6,14){\line(1,0){16}}

\put(10,3){\line(0,1){4}}
\put(12,3){\line(0,1){4}}
\put(16,3){\line(0,1){4}}
\put(19,3){\line(0,1){4}}

\put(10,9){\line(0,1){6}}
\put(12,9){\line(0,1){6}}
\put(16,9){\line(0,1){6}}
\put(19,9){\line(0,1){6}}

\put(8,4){\circle*{0.3}}
\put(10,6){\circle*{0.3}}
\put(12,6){\circle*{0.3}}

\put(8,2){\makebox(0,0){$\scriptstyle \infty$}}
\put(10,2){\makebox(0,0){$\scriptstyle 0$}}
\put(12,2){\makebox(0,0){$\scriptstyle 1$}}
\put(17.5,2){\makebox(0,0){$\scriptstyle \lambda_j\,\not=\,0,1,\infty$}}

\put(17,5){\makebox(1,1){$\cdots$}}
\put(17,12){\makebox(1,0){$\cdots$}}

\put(7,11.9){\makebox(0,1){$\vdots$}}
\put(14,11.9){\makebox(0,1){$\vdots$}}
\put(21,11.9){\makebox(0,1){$\vdots$}}

\put(1.5,11.9){\makebox(1,1){$\Xb$}}
\put(1.5,3.5){\makebox(1,1){$\bar{Y}$}}
\put(0,7.7){\makebox(1,1){$\bar{f}$}}
\put(24,3.5){\makebox(1,1){$\scriptstyle \Yb_0=\PP^1_k$}}
\put(2,10){\vector(0,-1){3.5}}

\end{picture}

\caption{\label{stablepic} The stable reduction of a three point cover}
\end{center}
\end{figure}
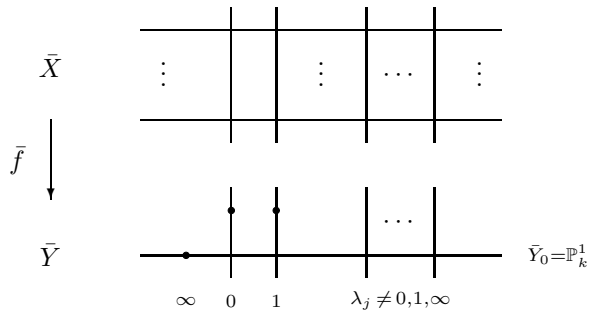

The following theorem is the first main result of \cite{bad}.

\begin{thm} \label{reductionthm}
  Suppose that $p$ strictly divides the order of $G$ and that $f$ has
  bad reduction. Then the following holds (compare with Figure
  \ref{stablepic}). 
  \begin{enumerate}
  \item The original component $\Yb_0$ is the only bad component.
    Every good component $\Yb_i$ intersects $\Yb_0$ in a unique point
    $\lambda_i\in\Yb_0$.
  \item The inertia group $I_0$ corresponding to the bad component
    $\Yb_0$ is cyclic of order $p$. The subcover $\Zb_0\to\Yb_0$ of
    $\fb_0$ (which is Galois with group $G_0/I_0$), is ramified at
    most in the points $\lambda_i$ (where $\Yb_0$ intersects a good
    component).
  \item
    For $i=1,\ldots,r$, the Galois cover $\fb_i:\Xb_i\to\Yb_i$ is
    wildly ramified at the point $\lambda_i$ and tamely ramified above
    $\Yb_i-\{\lambda_i\}$. If $\fb_i:\Xb_i\to\Yb_i$ is ramified at a
    point $\not=\lambda_i$, then this point is the specialization of
    one of the three branch points $0$, $1$, $\infty$ of the cover
    $f:X\to\PP^1_K$. 
  \end{enumerate}
\end{thm}

Part (ii) and (iii) of this theorem follow from part (i), by the
results of \cite{Raynaud98}. In fact, this implication is not
restricted to three point covers but holds for much more general
covers $f:X\to Y$. On the other hand, the truth of part (i) 
depends in an essential way on the assumption that $f$ is a three
point cover.

Under the additional assumption that all the ramification indices of
$f$ are prime to $p$, Theorem \ref{reductionthm} follows already from
the results of \cite{special}, via Raynaud's construction of the {\em
  auxiliary cover} (see the introduction of \cite{special}). 

Here is a brief outline of the proof of Theorem \ref{reductionthm}.
First, certain general results on the stable reduction of Galois
covers, proved in \cite{Raynaud98}, already impose severe restrictions
on the map $\fb$. For instance, it is shown that the good components
are precisely the {\em tails} of $\Yb$ (i.e.\ the leaves of the tree
of components of $\Yb$). Also, the Galois covers $\fb_i:\Xb_i\to\Yb_i$
(if $\Yb_i$ is good) and $\Zb_i\to\Yb_i$ (if $\Yb_i$ is bad) are
ramified at most at the points which are either singular points of the
curve $\Yb$ or specialization of a branch point of $f$.

In the next step one defines, for each bad component $\Yb_i$, a certain
differential form $\omega_i$ on the Galois cover $\Zb_i\to\Yb_i$. This
differential form satisfies some very special conditions, relative to the map
$\fb:\Yb\to\Xb$ and the action of $G$ on $\Yb$. For instance, $\omega_i$ is
either logarithmic (i.e.\ of the form ${\rm d}u/u$) or exact (i.e.\ of the
form ${\rm d}u$).  Furthermore, $\omega_i$ is an eigenvector under the action
of the Galois group $G_i/I_i$ of the cover $\Zb_i\to\Yb_i$, and its zeros and
poles are related to and determined by the ramification of the map
$\fb:\Yb\to\Xb$. These properties follow from the work of Henrio
\cite{YannickArbres}. Let us say for short that $\omega_i$ is {\em compatible}
with $\fb$.  Intuitively, $\omega_i$ encodes infinitesimal information about
the action of the inertia group $I_i$ on the stable model $\X$, in a
neighborhood of the component $\Xb_i$. Within the proof of Theorem
\ref{reductionthm}, the important point is that the existence of the
compatible differentials $\omega_i$ imposes further restrictions on the map
$\fb:\Xb\to\Yb$. In fact, these restrictions are strong enough to prove part
(i) of Theorem \ref{reductionthm}. For details, see \cite{bad}, \S 2.1.

\subsection{Special deformation data} \label{special}

By Theorem \ref{reductionthm} (i), the original component $\Yb_0=\PP^1_k$ is
the only bad component for the stable reduction of the three point cover $f$.
The proof of Theorem \ref{reductionthm} shows that there exists a differential
form $\omega_0$ on the Galois cover $\Zb_0\to\Yb_0$ which is compatible with
$\fb$, in the sense explained above. It is worthwhile to write down explicitly
what `compatibility' implies for the differential $\omega_0$. 

To simplify the exposition, we assume that the ramification indices of $f$ are
all divisible by $p$. If this is the case, then the branch points $0$, $1$ and
$\infty$ specialize to the original component $\Yb_0$. Since we identify
$\Yb_0$ with $\PP^1_k$, this means that the points
$\lambda_1,\ldots,\lambda_r$ where $\Yb_0$ intersects the good components
$\Yb_1,\ldots,\Yb_r$ are distinct from $0$, $1$, $\infty$.  By Theorem
\ref{reductionthm} (iii), the Galois covers $\fb_i:\Xb_i\to\Yb_i$ are then
\'etale over $\Yb_i-\{\lambda_i\}$.

Let $t$ denote the rational function on the original component $\Yb_0$
which identifies it with $\PP^1_k$. Compatability of $\omega_0$ with $\fb$
implies that 
\begin{equation} \label{omegaeq}
   \omega_0 \;=\; c\cdot\frac{z\,{\rm d}t}{t\,(t-1)},
\end{equation}
where $c\in k^\times$ is a constant and $z$ is a rational function on
$\Zb_0$ for which an equation of the form
\begin{equation} \label{zeq}
  z^{p-1} \;=\; \prod_i\, (t-\lambda_i)^{a_j}
\end{equation}
holds. Here the $a_i$ are integers $1<a_i<p$ such that $\sum_i
a_i=p-1$. These integers are determined by the (wild) ramification of the
Galois covers $\fb_i:\Xb_i\to\Yb_i$ at $\lambda_i$. 

Compatibility of $\omega_0$ with $\fb$ also implies that $\omega_0$ is
logarithmic, i.e.\ is of the form ${\rm d}u/u$, for some rational function $u$
on $\Zb_0$. Equivalently, $\omega_0$ is invariant under the Cartier operator.
The latter condition gives a finite list of equations satisfied by the
$t$-coordinates of the points $\lambda_i$ (depending on the numbers $a_i$).
One can show that these equations have only a finite number of solutions
$(\lambda_i)$, see \cite{cotang}, Theorem 5.14. In other words, the existence
of the differential form $\omega_0$ determines the position of the points
$\lambda_i$, up to a finite number of possibilities.

The pair $(\Zb_0,\omega_0)$ is called a {\em special deformation datum}. Given
a special deformation datum $(\Zb_0,\omega_0)$, the branch points $\lambda_i$
of the cover $\Zb_0\to\Yb_0$ are called the {\em supersingular points}. A
justification for this name, in form of a well known example, will be given in
\S \ref{modular}. By the result mentioned in the preceeding paragraph, a
special deformation datum is {\em rigid}, i.e.\ is an object with
$0$-dimensional moduli.  This is no surprise, as special deformation data
arise from three point covers, which are rigid objects themselves.  We point
this out, because this sort of rigidity is the (somewhat hidden) principle
underlying all results discussed in this note which are particular for three
point covers.  In the following section, we will interpret the existence of
$(\Zb_0,\omega_0)$ as a liftability condition for the map $\fb:\Xb\to\Yb$.

\section{Lifting} \label{lifting}

\subsection{Special $G$-maps} \label{lift}

The stable reduction of a three point Galois cover $f:X\to\PP^1_K$ is,
by definition, a finite map $\fb:\Xb\to\Yb$ between semistable curves
over the residue field $k$, together with an embedding
$G\inj\Aut(\Xb/\Yb)$. In the case of bad reduction, the curves $\Yb$
and $\Xb$ are singular, and the map $\fb$ is inseparable over some of
the components of $\Yb$. This suggests the following question. Given a
map $\fb:\Xb\to\Yb$ of the sort we have just described, together with
an embedding $G\inj\Aut(\Xb/\Yb)$, does it occur as the stable
reduction of a three point Galois cover $f:X\to\PP^1_K$, for some
finite extension $K/K_0$? If this is the case, then we say that
$f:X\to\PP^1_K$ is a {\em lift} of $\fb:\Xb\to\Yb$.

Theorem \ref{reductionthm} and its proof give a list of necessary conditions
on $\fb$ for the existence of a lift (at least under the extra condition that
$p$ strictly divides the order of $G$).  These conditions lead naturally to
the notion of a {\em special $G$-map}. See \cite{bad}, \S 2.2 for a precise
definition. To give the general idea, it suffices to say that a special
$G$-map is a finite map $\fb:\Xb\to\Yb$ between semistable curves, together
with an embedding $G\inj\Aut(\Xb/\Yb)$, which admits a compatible special
deformation datum $(\Zb_0,\omega_0)$. One can show that special $G$-maps are
rigid in the sense we used at the end of \S \ref{special}. Moreover, one has
the following lifting result, proved in \cite{bad}, \S 4.

\begin{thm} \label{liftthm} 
  Let $\fb:\Xb\to\Yb$ be a special $G$-map over $k$. Then the
  following holds.
  \begin{enumerate}
  \item
    There exists a three point cover $f:X\to\PP^1$ lifting $\fb$.
  \item
    Every lift $f$ of $\fb$ can be defined over a finite extension $K/K_0$
    which is at most tamely ramified.
  \end{enumerate}
\end{thm}

The corresponding result proved in \cite{bad}, \S 4, is somewhat
stronger. It determines the set of isomorphism classes of all lifts of
$\fb$, together with the action of $\Gal(\Kb_0/K_0)$, in terms of
certain invariants of $\fb$ (these invariants are essentially the
numbers $a_i$ appearing in \eqref{zeq}). This more precise result
gives an upper bound for the degree of the minimal extension $K/K_0$
over which every lift of $\fb$ can be defined.

Theorem \ref{tamethm} follows easily from Theorem \ref{reductionthm} and
Theorem \ref{liftthm} (in a way similar to how Beckmann's Theorem follows from
Grothendieck's theory of tame covers, see \S \ref{stable}).

Part (i) of Theorem \ref{liftthm}, i.e.\ the mere existence of a lift,
follows already from the results of \cite{special}. Part (ii) is more
difficult. The technical heart of the proof is a study of the
deformation theory of a certain curve with an action of a finite group
scheme, which is associated to a special deformation datum. A detailed
exposition of this deformation theory can be found in \cite{cotang}.
An overview of the proof of Theorem \ref{liftthm} will be given in \S
\ref{outline} below.

\subsection{The supersingular disks} 
   \label{supersingular}

Let $f:Y\to\PP^1_K$ be a three point cover, defined over a finite
extension $K/K_0$, with bad reduction. Let $\fb:\Xb\to\Yb$ be the
stable reduction of $f$. We assume that the conclusion of Theorem
\ref{reductionthm} holds (it holds, for instance, if $p^2\nmid |G|$).
Let $Y\an$ denote the rigid analytic $K$-space associated to
$Y=\PP^1_K$. The $\OO_K$-model $\Y$ of $Y$ yields a specialization map
${\rm sp}_{\Y}:Y\an\to\Yb$. For $i=1,\ldots,r$, the good component
$\Yb_i$ gives rise to a rigid analytic subset
\[
    D_i \;:=\; {\rm sp}_{\Y}^{-1}(\Yb_i-\{\lambda_i\}) \;\subset\; Y\an.
\]
As a rigid $K$-space, $D_i$ is a closed unit disk, i.e.\ is isomorphic
to the affinoid ${\rm Spm}\,K\{\{T\}\}$. See e.g.\ \cite{HenrioDisks}. 

An important step in the proof of Part (ii) of Theorem \ref{liftthm}
is to show that the disks $D_i$ depend only on the reduction $\fb$,
but not on the lift $f$ of $\fb$. For simplicity, we shall again
assume that all the ramification indices of the three point cover
$f:X\to\PP^1_K$ are divisible by $p$, see \S \ref{special}. Then the
special deformation datum $(\Zb_0,\omega_0)$ associated to the
reduction $\fb:\Xb\to\Yb$ is essentially determined by points
$\lambda_1,\ldots,\lambda_r\in\PP^1_k-\{0,1,\infty\}$ and integers
$a_1,\ldots,a_r$ with $1<a_i<p$ and $\sum_ia_i=p-1$. We consider
$\lambda_i$ as an element of $k-\{0,1\}$ and let
$\tilde{\lambda}_i\in\OO_{K_0}$ be a lift of $\lambda_i$. By
definition, the closed disk $D_i$ is contained in the open unit disk
\[
    D_i' \;:=\; \{t\in\OO_K\,\mid\, |t-\tilde{\lambda}_i|_K<1\,\}
       \subset Y\an. 
\]
With this notation, we have the following result, see \cite{bad},
Proposition 4.3 and the remark following the proof of Theorem 3.8.

\begin{prop} \label{diskprop} 
  We have 
  \[
       D_i \;=\; \{\,t\in D_i'  \;\mid\; |t-\tilde{\lambda}_i|_K 
           \,\leq\, |p|_K^\frac{p}{p-1+a_i}\,\}.
  \]
  In particular, the disks $D_i$ depend only on the special
  deformation datum $(\Zb_0,\omega_0)$.
\end{prop}

We call the open disks $D_i'$ the {\em supersingular disks} associated
to the special deformation datum $(\Zb_0,\omega_0)$. This is in
correspondence with naming the points $\lambda_i\in\Yb_0$ the
supersingular points, see \S \ref{special}. The closed subdisk
$D_i\subset D_i'$ is called the {\em too supersingular disk}, a term
which is also borrowed from the theory of moduli of elliptic curves,
see \S \ref{modular}.

\subsection{The auxiliary cover} \label{aux}

We continue with the notation and assumptions of the preceeding
subsection.  Recall that we have chosen in \S \ref{badreduction} a
component $\Xb_i\subset\Xb$ above the component $\Yb_i\subset\Yb$. The
stabilizer of $\Xb_i$ is the subgroup $G_i\subset G$ and the inertia
subgroup of $\Xb_i$ is a normal subgroup $I_i\lhd G_i$. By the
conclusion of Theorem \ref{reductionthm}, we have $|I_0|=p$ and
$|I_i|=1$ for $i=1,\ldots,r$. Let $\Xb\hor\subset\Xb$ denote the union
of all components mapping to the original component $\Yb_0$. Then we
have
\[
   f^{-1}(D_i) \;=\; \Ind_{G_i}^G(E_i), \quad\text{\rm with}\;\;
      E_i \;=\; {\rm sp}_{\X}^{-1}(\Xb_i-\Xb\hor) \;\subset\; X\an.
\]
The map $E_i\to D_i$ is a finite \'etale Galois cover between smooth
affinoid $K$-spaces whose reduction is equal to the restriction of the
\'etale Galois cover
$\fb_i^{-1}(\Yb_i-\{\lambda_i\})\to\Yb_i-\{\lambda_i\}$. This
determines $E_i\to D_i$ uniquely, up to isomorphism, because lifting
of \'etale morphisms is unique.

Let $U_0:=Y\an-(\cup_i D_i)$ denote the complement of the disks
$D_i$. Then we have
\[
   f^{-1}(U_0) \;=\; \Ind_{G_0}^G(V_0), \quad\text{\rm with}\;\;
   V_0 \;:=\; {\rm sp}_{\X}^{-1}(\Xb\hor) \;\subset\; X\an.
\]
The map $V_0\to U_0$ is a finite Galois cover between smooth (non
quasicompact) rigid $K$-spaces, \'etale outside the subset
$\{0,1,\infty\}\subset U_0$. It can be shown that there exists a
$G_0$-Galois cover $f\aux:X\aux\to Y=\PP^1_K$ such that
$V_0=(f\aux)^{-1}(U_0)$. Such a cover $f\aux$ is ramified at $0$, $1$,
$\infty$ and, for each $i=1,\ldots,r$, at one point $y_i\in D_i'$.
Moreover, the cover $f\aux$ is uniquely determined by the choice of
the points $y_i$. It is called the {\em auxiliary cover} associated to
$f$ and the points $(y_i)$, see \cite{Raynaud98}, \cite{special} and
\cite{bad}, \S 4.1.3.  Let $\partial D_i$ denote the boundary of the
disk $D_i$. By construction of $f\aux$, there exists a $G$-equivariant
isomorphism
\begin{equation} \label{patchingeq}
  \varphi_i:\Ind_{G_i}^G\,(E_i\times_{D_i}\partial D_i) \;\cong\;
   \Ind_{G_0}^G\, (f\aux)^{-1}(\partial D_i),
\end{equation}
compatible with the natural map to $\partial D_i$.

\subsection{The proof of Theorem \ref{liftthm}} \label{outline}

We will now give a brief outline of the proof of Theorem \ref{liftthm}.
Suppose that $\fb:\Xb\to\Yb$ is a special $G$-map. We want to construct all
three point covers $f:X\to Y=\PP^1_K$ lifting $\fb$. For the moment, we let
$K$ be any sufficiently large finite extension of $K_0$. At the end of our
argument, we will reason that it suffices to take for $K$ a certain tame
extension of $K_0$, which is explicitly determined by $\fb$.

We divide the proof into three steps.  The first step consists in
constructing a Galois cover $f\aux:X\aux\to Y$ which can play the role
of the auxiliary cover associated to any lift $f$ of $\fb$. In fact,
one shows that all good candidates $f\aux$ live in a continuous family
which depends only on the special deformation datum $(\Zb_0,\omega_0)$
associated to $\fb$. The individual members of this family depend on
the choice of the extra branch points $y_i\in D_i'$ (with $D_i'$ as
above). See \cite{special}, \S 3.2, and \cite{bad}, \S 3.

Let $D_i\subset D_i'$ be a smaller closed disk, defined as in
Proposition \ref{diskprop} (the numbers $a_i$ used in the
definition of $D_i$ are determined by the special deformation datum
$(\Zb_0,\omega_0)$).  Let $E_i\to D_i$ be the \'etale $G_i$-Galois
cover lifting $\fb_i^{-1}(\Yb_i-\{\lambda_i\})\to\Yb_i-\{\lambda_i\}$.
For any choice of points $y_i\in D_i'$, we obtain a cover
$f\aux:X\aux\to Y$, which is a candidate for the auxiliary cover. 
In this situation, a tuple $(\varphi_i)$ of isomorphisms as in
\eqref{patchingeq} is called a {\em patching datum}. 

The second step of the proof consists in showing that there exists a
patching datum $(\varphi_i)$ if and only if the points $y_i$ lie in
the smaller disk $D_i$. The sufficiency of the condition $y_i\in D_i$
can be shown using the same arguments as in \cite{special}, \S 3.4.
The necessity of this condition -- which is equivalent to Proposition
\ref{diskprop} above -- lies somewhat deeper.  See \cite{bad}, \S 3,
in particular Theorem 3.8.  In this step one uses the deformation
theory developed in \cite{cotang}.

The third and final step uses rigid (or formal) patching.  For any
choice of $y_i\in D_i$, let $f\aux:X\aux\to Y=\PP^1_K$ be the
associated auxiliary cover, and set $V_0:=(f\aux)^{-1}(U_0)$. By the
second step, we have a patching datum $(\varphi_i)$. The proof of the
claim in step two shows moreover that the cover $V_0\to U_0$ depends
only on the special deformation datum, but not on the choice of
$y_i\in D_i$. Using rigid patching, one easily constructs a $G$-Galois
cover $f:X\to Y$ such that $f^{-1}(D_i)=\Ind_{G_i}^G(E_i)$,
$f^{-1}(U_0)=\Ind_{G_0}^G(V_0)$ and such that the patching datum
$(\varphi_i)$ is induced by the identity on $X$.  Essentially by
construction, $f$ is a three point cover lifting the special $G$-map
$\fb$. This proves Part (i) of Theorem \ref{liftthm}.

It is not hard to see that all lifts of $\fb$ arise in the way we have just
described. More precisely, the set of isomorphism classes of lifts of $f$ is
in bijection with the set of patching data. Therefore, to finish the proof of
Theorem \ref{liftthm} it suffices to show that the above construction works
over a tame extension $K/K_0$. Actually, the construction of the covers
$E_i\to D_i$ and of the auxiliary cover $f\aux:X\aux\to Y$ can be done over
$K_0$ (set $y_i:=\tilde{\lambda}_i$). A direct analysis shows that patching
data $(\varphi_i)$ exist if one takes for $K$ the (unique) tame extension of
degree $(p-1)\cdot{\rm lcm}_i(p-1+a_i)$.  This concludes the proof of Theorem
\ref{liftthm}.

\section{Modular curves} \label{modular}

\subsection{Modular curves as three point covers}

Let $\HH$ denote the upper half plane and $\HH^*$ the union of $\HH$
with the set $\PP^1(\QQ)$. The group $\SL_2(\ZZ)$ acts on $\HH$ and
$\HH^*$ in a standard way. Moreover, for every subgroup
$\Gamma\subset\SL_2(\ZZ)$ of finite index, the quotient
$X_\Gamma:=\HH^*/\Gamma$ carries a natural structure of a compact
Riemann surface, and therefore also of a smooth projective curve over
$\CC$. The classical $j$-function identifies the quotient of $\HH^*$
by $\SL_2(\ZZ)$ with the projective line. So for each finite index
subgroup $\Gamma$ we obtain a finite cover of compact Riemann surfaces
(or smooth projective curves over $\CC$)
\[
      f_\Gamma:\,X_\Gamma \;\To\; \PP^1.
\]
This map is unramified away from the three points $0,1728,\infty$. In
other words, $f_\Gamma$ is a three point cover. 

For an integer $N$, let 
\[
       \Gamma(N) \;:=\; \{\; A\in\SL_2(\ZZ) \;\mid\;
           A \,\equiv\, I_2 \pmod{N} \;\}.
\]
A {\em congruence subgroup} is a subgroup $\Gamma\subset\SL_2(\ZZ)$
which contains $\Gamma(N)$, for some $N$.  The corresponding curve
$X_\Gamma$ is called a {\em modular curve}.  The standard examples for
congruence subgroups are $\Gamma(N)$, $\Gamma_0(N)$ and $\Gamma_1(N)$.
The corresponding modular curves are usually denoted by $X(N)$,
$X_0(N)$ and $X_1(N)$.

\subsection{The modular curve $X(p)$}

Let us fix an odd prime number $p$. The three point cover $g:X(p)\to
X(1)=\PP^1$ is Galois, with Galois group $G=\PSL_2(p)$. Its ramification
index at the branch point $\infty$ (resp.\ at $0$, resp.\ at 
$1728$) is equal to $p$ (resp.\ $3$, resp.\ $2$). Note that the order
of $G$ is equal to $p(p^2-1)/2$. 

For technical reasons it is easier to discuss a variant of $g$, namely
the three point cover
\[
     f:\,X(2p) \;\to\; X(2) \,=\, \PP^1.
\]
(We identify the modular curve $X(2)$ with $\PP^1$ by means of the classical
$\lambda$-function.) The cover $f$ is Galois with Galois group $G=\PSL_2(p)$.
The ramification index at each of the three branch points $0$, $1$, $\infty$
is equal to $p$. There is an equivariant action of $\SL_2(2)\cong S_3$ on the
source and the target of $f$. The cover $g$ is obtained by taking the quotient
under this action.

The next proposition follows easily from the results of
\cite{DeligneRapoport}.

\begin{prop} \label{modularprop}
  The covers $f$ and $g$ can be defined over $\QQ$. In particular,
  $\QQ$ is their field of moduli.  Furthermore, $g$ (resp.\ $f$)
  has good reduction at the prime $l$ (in the sense of Definition
  \ref{stabdef}) for $l\not=2,3,p$ (resp.\ for $l\not=p$).
\end{prop}

The fact that $f$ and $g$ can be defined over $\QQ$ follows
also from the {\em rigidity criterion} used in inverse Galois theory.
See \cite{MalleMatzat} or \cite{Voelklein}. The $\QQ$-models of $f$
and $g$ are not unique. However, both covers have a unique model over
the field $\QQ(\sqrt{p^*})$ which is Galois.

Note that the second statement of Proposition \ref{modularprop}
confirms (but is not implied by) the good reduction criterion
discussed in \S \ref{good}. For $l=2,3$, the curve $X(p)$ has good
reduction as well.  However, since in this case $l$ divides one of the
ramification indices, the cover $g$ has bad reduction -- at least in
the sense of Definition \ref{stabdef}.

Since $p$ exactly divides the order of $G$, the results discussed in
\S \ref{reduction} and \S \ref{lifting} can be used to study the
stable reduction of $f$ and $g$ at the prime $p$. This in done in
detail in \cite{modular}. We shall present some of the main results of
\cite{modular}, as an illustration for the results discussed earlier. 

Let $K_0$ be the completion of the maximal unramified extension of
$\QQ_p$. From now on, we consider the three points covers $f$ and
$g$ as defined over $K_0$. We remark that there are many different
models of these covers over $K_0$. However, the stable model, which
exists over a finite extension of $K_0$, is unique. 

For simplicity, we discuss only the stable reduction of the cover $f$
in detail. Let $K/K_0$ be the minimal extension over which $f$ has
semistable reduction, and let $\fb:\Xb\to\Yb$ denote the stable
reduction of $f$.  We will freely use the notation introduced in \S
\ref{badreduction} and \S \ref{special}.  Since $f$ has bad reduction
and $p$ strictly divides the order of $G$, the conclusion of Theorem
\ref{reductionthm} holds.  In particular, $\fb$ gives rise to a
special deformation datum $(\Zb_0,\omega_0)$ and to Galois covers
$\fb_i:\Xb_i\to\Yb_i$. Since all the ramification indices of the cover
$f$ are equal to $p$, none of the supersingular points
$\lambda_i\in\Yb_0=\PP^1_k$ equal $0$, $1$ or $\infty$.
Moreover, the Galois covers $\fb_i$ are \'etale over
$\Yb_i-\{\lambda_i\}$.

\begin{thm}\label{modularthm} \ 
  \begin{enumerate}
  \item The field $K$ is the (unique) tame extension of $K_0$ of degree
    $(p^2-1)/2$.
  \item There are $r=(p-1)/2$ supersingular points $\lambda_i$; they are
    precisely the roots of the {\em Hasse polynomial}
    \[
          \Phi(t) \;=\; \sum_{j=0}^{r}\; \binom{r}{j}^2\; t^j.
    \]
    The Galois cover $\Zb_0\to\Yb_0$ is the cyclic cover of degree $r$
    given by the equation $z^r \;=\; \Phi(t)$. Furthermore, we have
    $\omega_0=z\,t^{-1}(t-1)^{-1}{\rm d}t$.
  \item
    For $i=1,\ldots,r$, the curve $\Xb_i$ is given by the equation
    $y^{(p+1)/2}\,=\,x^p\,-\,x$.
    An element of $G=\PSL_2(p)$, represented by a matrix
    $A=\left(\begin{smallmatrix} a & b \\ c & d 
     \end{smallmatrix}\right)\in\SL_2(p)$, acts on $\Xb_i$ as follows:
    \[
        A(x) \;=\; \frac{ax+b}{cx+d}, \qquad 
        A(y) \;=\; \frac{y}{(cx+d)^2}.
    \]
  \end{enumerate}
\end{thm}
         
For the proof of this theorem, see \cite{modular}. The main idea is this. One
constructs a special $G$-map $\fb:\Xb\to\Yb$, which satisfies the conclusion
of Theorem \ref{modularthm} (ii) and (iii).  By Theorem \ref{liftthm}, it
lifts to a three point cover $f':X'\to\PP^1$, defined over the tame extension
of $K_0$ of degree $(p^2-1)/2$. By the rigidity criterion of inverse Galois
theory, $f'$ has to be isomorphic to $f$.  Whence the theorem.

One can prove a similar theorem about the stable reduction of $g:X(p)\to
X(1)=\PP^1$. If $p\equiv 1\pmod{12}$ then the statements are almost identical,
except that the Hasse polynomial has to be replaced by another polynomial (for
which there is an explicit formula, similar to the expression for $\Phi$). If
$p\not\equiv 1\pmod{12}$ then some of the supersingular points $\lambda_i$ are
equal to either $0$ or $1728$, and the corresponding Galois cover
$\fb_i:\Xb_i\to\Yb_i$ is not the one from Theorem \ref{modularthm} (iii).

The modular curves $X_0(p)$, $X_1(p)$ and $X_0(p^2)$ are all quotients of the
curve $X(p)$. Using the results of \cite{modular} on the stable reduction of
$X(p)$, one can determine the stable reduction of all these quotients,
reproving results of Deligne--Rapoport \cite{DeligneRapoport} and Edixhoven
\cite{Edixhoven90}. Somewhat surprisingly, this new proof does not use the
interpretation of modular curves as moduli spaces for elliptic curves with
level structure. 

However, the `modular' interpretation of modular curve justifies the use of
the term `supersingular' in \S \ref{special} and \S \ref{supersingular}. In
fact, the supersingular points $\lambda_i$, which are the roots of the Hasse
polynomial $\Phi$, are exactly the values $t\in k$ for which the Legendre
elliptic curve $E_t$ with equation $y^2=x(x-1)(x-t)$ is supersingular.
Similarly, a point $t\in\PP^1(K)$ lies in one of the open disks $D_i'$ if and
only if the elliptic curve $E_t$ has supersingular reduction. The modular
interpretation of the smaller closed disks $D_i\subset D_i'$ is somewhat less
known. However, it can be shown (see e.g.\ \cite{Coleman}) that a point
$t\in\PP^1(K)$ lies in one of the disks $D_i$ if and only if the elliptic
curve $E_t$ is {\em too supersingular}, which means that $E_t[p]$ has no
canonical subgroup.

\end{document}